\documentclass{amsproc}

\usepackage{amsmath, amssymb, amsthm, xcolor,enumerate}

\usepackage{hyperref}
\usepackage{xypic}
\usepackage{graphicx}

\newtheorem{theorem}{Theorem}[section]

\theoremstyle{definition}
\newtheorem{definition}[theorem]{Definition}

\theoremstyle{remark}

\numberwithin{equation}{section}

\renewcommand{\H}{\mathcal{H}}
\newcommand{\A}{\mathcal{A}}
\newcommand{\B}{\mathcal{B}}
\newcommand{\tr}{\mathrm{Tr}}

\begin{document}

\title{$\mathrm{C}^*$-algebras and their nuclear dimension}

\author{Jorge Castillejos}

\address{Department of Mathematics, KU Leuven, 3001 Leuven, Belgium.}
\curraddr{Institute of Mathematics, Polish Academy of Sciences, ul. {\'S}niadeckich 8, 00-656 Warszawa, Poland}

\email{jcastillejoslopez@impan.pl}

\thanks{The author is supported by long term structural funding - Methusalem grant of the Flemish Government.}


\begin{abstract}
We review the notion of nuclear dimension for $\mathrm{C}^*$-algebras introduced by Winter and Zacharias. We explain why it is a non-commutative version of topological dimension. After presenting several examples, we give a brief overview of the state of the art.
\end{abstract}

\maketitle

\section*{Introduction}

Typically, \textit{$\mathrm{C}^*$-algebras} are defined as complex Banach algebras with an involution that satisfy the $\mathrm{C}^*$-identity. Precisely, a complex algebra $\A$ with an involution is a $\mathrm{C}^*$-algebra if it has a complete norm that makes the operations continuous and satisfies the $\mathrm{C}^*$-\textit{identity},
\begin{equation}
\| a^* a \| = \| a \|^2, \qquad a \in \A. \notag
\end{equation}
Thanks to the Gelfand-Naimark theorem, we can also define $\mathrm{C}^*$-algebras as subalgebras of the algebra of bounded operators on a Hilbert space \cite{GN43}. 
This immediately leads us to the theory of \textit{operator algebras}. 

Operator algebras are subalgebras of the algebra of bounded operators on a Hilbert space. In particular, $\mathrm{C}^*$-algebras are those that are closed and self-adjoint. 
Motivated by the fact that unital commutative $\mathrm{C}^*$-algebras are always isomorphic to the algebra of complex-valued continuous functions on a compact Hausdorff space, we view general $\mathrm{C}^*$-algebras as non-commutative topological spaces.
This area of mathematics originated from the groundbreaking work of von Neumann on the mathematical formulation of quantum mechanics \cite{vNn55}, and nowadays has interesting connections with other areas of mathematics such as group theory, dynamical systems, geometry, logic, and quantum information. 

Many endeavours have been directed towards  understanding the structure of $\mathrm{C}^*$-algebras. In particular, the well-behaved class of \textit{nuclear} $\mathrm{C}^*$-algebras has been deeply studied in the past \cite{Tak64,Lan73,CE78,Ki77}. 
Roughly speaking, a $\mathrm{C}^*$-algebra is nuclear if it can be approximated by matrix algebras
and nuclearity can also be considered as the $\mathrm{C}^*$-version of \textit{amenability} for groups. 

By viewing $\mathrm{C}^*$-algebras as non-commutative topological spaces, Winter introduced several notions of non-commutative topological dimension for nuclear $\mathrm{C}^*$-algebras using \v{C}ech's \textit{covering dimension} as a model \cite{Wi03, cov-dim2}. These notions were later refined by Winter and Kirchberg \cite{Kir-Win}, and by Winter and Zacharias \cite{WZ10}. 
Roughly speaking, the covering dimension of a topological space is obtained by colouring finite open covers of the space and finding the minimum number of colours that can be used in such a way that open sets with the same colour do not intersect each other; the dimension is equal to the minimum number of colours minus one. For instance, the covering dimension of the interval $[0,1]$ is equal to $1$ because we need at least $2$ colours. 

\begin{figure}[ht!]
	\centering
	\includegraphics[width=120mm]{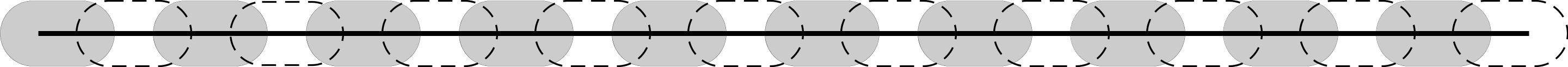}
\end{figure}

By Gelfand duality, unital commutative $\mathrm{C}^*$-algebras are always equal to the algebra of complex-valued continuous functions on some compact Hausdorff space $X$. The covering dimension of the underlying space $X$ is reflected by certain finite dimensional approximations of the algebra of continuous functions. These approximations are used to model several versions of non-commutative covering dimension. After several iterations, the definitive notion appeared as \textit{nuclear dimension}. This revolutionary notion became fundamental in the understanding of the fine structure of nuclear $\mathrm{C}^*$-algebras and played a major role in the \textit{classification programme} of simple separable unital nuclear $\mathrm{C}^*$-algebras. We refer to \cite{WiSurvey} for an overview.

The classification programme of nuclear $\mathrm{C}^*$-algebras was initiated by Elliott in his seminal work \cite{Ell76} where he classified all unital \textit{approximately finite dimensional} $\mathrm{C}^*$-algebras via \textit{ordered K-theory}. He managed to extend his classification theorem to larger classes of $\mathrm{C}^*$-algebras \cite{Ell93}. Based on this, he was prompted to conjecture that we can classify simple separable nuclear and unital $\mathrm{C}^*$-algebras with K-theory and \textit{traces} \cite{elliott-classification}. For technical reasons, we must also assume that the $\mathrm{C}^*$-algebras satisfy the \textit{Universal Coefficient Theorem} (UCT). This technical condition can roughly be considered as being ``very weak homotopy equivalent'' to a commutative $\mathrm{C}^*$-algebra. It is unknown if there is a separable nuclear $\mathrm{C}^*$-algebra that does not satisfy the UCT.

To the surprise of many, this conjecture was verified for large classes of $\mathrm{C}^*$-algebras \cite{DE02,EGLP96,EG96,EV93,AH-classification,Lin03,Lin04,Lin11,Wi14,Wi16}, including the spectacular Kirchberg-Phillips classification of Kirchberg algebras that satisfy the UCT \cite{Phi00}.
However, counterexamples to this conjecture were found by R\o rdam and Toms \cite{Ror03,Toms}. It became clearer then that some regularity condition must be added to the conjecture \cite{regularity-properties}. 

The missing ingredient in Elliott's conjecture was precisely finite nuclear dimension. Indeed, the aforementioned counterexamples to the conjecture have infinite nuclear dimension!
After several decades of work, one of the major achievements in $\mathrm{C}^*$-algebra theory was completed: the classification of simple separable unital $\mathrm{C}^*$-algebras with finite nuclear dimension that satisfy the UCT \cite{GLN,EGLN,TWW}.

In this document, we will review Winter and Zacharias' notion of nuclear dimension. We will make a detailed analysis of the commutative case as a motivation for the general definition. Then, we provide several examples of $\mathrm{C}^*$-algebras and discuss their nuclear dimension. 

Let us finish this introduction by explaining the structure of this document. In Section 1 we review basic notions of $\mathrm{C}^*$-algebras. In Section 2 we discuss covering dimension and its non-commutative version, nuclear dimension. After this, we introduce several examples of $\mathrm{C}^*$-algebras and their nuclear dimension: approximately finite dimensional $\mathrm{C}^*$-algebras (Section 3), group $\mathrm{C}^*$-algebras (Section 4), $\mathrm{C}^*$-algebras associated to dynamical systems (Section 5) and Cuntz algebras (Section 6). In Section 7, we explain a trichotomy for nuclear dimension in the simple separable case. We finish this survey by proving that the nuclear dimension of $C(X)$ is equal to the covering dimension of the compact space $X$.

\subsection*{Acknowledgements} 
The author would like to thank Sam Evington, Jamie Gabe, Adrian Gonzalez-Perez and Stuart White
for their comments on an earlier version of this document. The author also thanks the anonymous referee for useful suggestions on an earlier version of this paper.

\section{Basics of $\mathrm{C}^*$-algebras}

Let us briefly introduce $\mathrm{C}^*$-algebras. A more comprehensive introduction can be found at \cite{Murphy}. 
Given a complex Hilbert space $\H$, a linear operator $T: \H \to \H$ is \emph{bounded}
if there exists a positive constant $M$ such that
\begin{equation}
\| T(h) \| \leq M \| h \|, \qquad \qquad h \in \H.
\end{equation}
It is well known that a linear operator is continuous if and only if it is bounded. Let us denote by $B(\H)$ to the set of bounded operators on $\H$. This set can be endowed with the structure of a complete normed algebra. Indeed, the sum and scalar multiplication are defined pointwise, the product is given by composition and the operator norm is defined as
\begin{equation}
\| T \| := \sup_{\| h \|=1} \|T(h)\|.
\end{equation}

It is an standard fact that $B(\H)$ is complete with respect to this norm. But $B(\H)$ has another interesting operation, a natural involution that is given by the adjoint. Given an operator $T \in B(\H)$, there exists a unique operator $T^* \in B(\H)$ that satisfies 
\begin{equation}
\langle T(h) , g \rangle = \langle h, T^* (g) \rangle, \qquad g,h \in \H.
\end{equation}

The operator $T^*$ is called the \textit{adjoint} of $T$ and the map $T \mapsto T^*$ defines an \textit{involution}, i.e.\ a conjugate linear map such that $(T^*)^* = T$ and $(ST)^* = T^* S^*$.

Another fundamental property of $B(\H)$ is the $\mathrm{C}^*$-\textit{identity},
\begin{equation}
\| T \|^2 = \| T^* T \| , \qquad \qquad T \in B(\H).
\end{equation} 
The importance of this identity relies on the fact that it forces the norm to be determined by the algebraic structure of $B(\H)$. Indeed, the norm of $T^*T$ depends solely on its \textit{spectrum} $\sigma(T^*T)$, 
the subset of all scalars $\lambda \in \mathbb{C}$ such that $T^*T - \lambda 1_{B(\H)}$ is not an invertible operator.

\begin{definition}
	A subalgebra $\A$ of $B(\H)$ is a \textit{$\mathrm{C}^*$-algebra} if it is closed with respect to
	the operator norm and closed under taking adjoints.
	
	\end{definition}

As already mentioned before, this is not the typical definition, but thanks to the Gelfand-Neimark theorem \cite[Section 3.4.1]{Murphy}, our definition agrees with the typical one. 
The $\mathrm{C}^*$-identity has strong consequences in the structure of $\mathrm{C}^*$-algebras, for instance the norm on a $\mathrm{C}^*$-algebra is unique and any $^*$-\textit{homomorphism}\footnote{A $^*$-\textit{homomorphism} $\varphi: \A \to \B$ is a  linear map such that $\varphi(ab)=\varphi(a)\varphi(b)$ and $\varphi(a^*)= \varphi(a)^*$ for al $a,b \in \A.$} $\varphi: \A \to \B$ is contractive and hence continuous.

An easy example of a $\mathrm{C}^*$-algebra is  $\mathbb{C}$, the complex numbers, where involution is simply given by complex conjugation. 
More generally, it follows from the definition that $B(\H)$ is a $\mathrm{C}^*$-algebra. 
When the Hilbert space is of finite dimension $n$, the algebra of bounded operators is simply 
$M_n(\mathbb{C})$, the algebra of $n \times n$ matrices with entries in $\mathbb{C}$, endowed with the operator norm, matrix multiplication as product and adjoint given by $(a_{ij})^*_{i,j} = (\overline{a_{ji}})_{i,j}$. In fact, any finite dimensional $\mathrm{C}^*$-algebra is always a direct sum of finitely many matrix algebras endowed with the operator norm. Notice that it is crucial to endow the operator norm.
For instance, $\mathbb{C}^2$ endowed with the norm $\|(z_1, z_2) \|_1 = |z_1|+|z_2|$ is not a $\mathrm{C}^*$-algebra but it is if instead we endow it with the supremum norm $\|(z_1, z_2) \|_\infty = \max \{|z_1|,|z_2|\}$.

\subsection{Commutative $\mathrm{C}^*$-algebras}

Given a compact Hausdorff space $X$, we can always construct a $\mathrm{C}^*$-algebra by considering
\begin{equation}
C(X) : = \{ f: X \to \mathbb{C} \mid f \text{\, is continuous} \}.
\end{equation} 
Operations are given by pointwise addition, pointwise multiplication, involution as $f^*(x)=\overline{f(x)}$
and norm given by 
\begin{equation}
\| f \|_{\infty} = \sup_{x \in X} |f(x)|.
\end{equation}
The algebra $C(X)$ acts on the Hilbert space $L^2(X, \mu)$ (where $\mu$ is a Borel measure) by multiplication.

A $\mathrm{C}^*$-algebra $\A$ is \textit{unital} if it has a multiplicative unit, denoted by $1_{\A}$.
Observe that the algebra $C(X)$ is unital, where the unit $1_{C(X)}$ is the constant function $1$. 
However, not every $\mathrm{C}^*$-algebra is unital. For instance, if $X$ is a locally compact Hausdorff space which is not compact, the corresponding algebra $C(X)$ is not normed \cite[Proposition 2.1.14]{Da00} so we have to consider some subalgebra instead. 

A continuous function $f:X \to \mathbb{C}$ \textit{vanishes at infinity} if for all $\epsilon>0$ there is a compact subset $K \subseteq X$ for which $|f(x)|\leq \epsilon$ if $x \in X\setminus K$. Set
\begin{equation}
C_0(X):= \{ f \in C(X) \mid f \text{ vanishes at infinity}\}.
\end{equation}
In this subalgebra, $\| \cdot \|_\infty$ is a norm and $C_0(X)$ is indeed a non-unital commutative $\mathrm{C}^*$-algebra since the constant function $1$ does not vanishes at infinity. Observe that if $X$ is compact, then $C(X)= C_0(X)$. 

We now state Gelfand's theorem that characterises commutative $\mathrm{C}^*$-algebras. 
We refer to \cite[Theorem 2.1.10]{Murphy} for the details of its proof.
Two $\mathrm{C}^*$-algebras $\A$ and $\B$ are \textit{isomorphic}, $\A \cong \B$, if there is an $^*$-isomorphism $\varphi: \A \to \B$.

\begin{theorem}[Gelfand] 
	If $\A$ is a commutative $\mathrm{C}^*$-algebra, then $\A \cong C_0(X)$ for some locally compact Hausdorff space $X$. In particular, $\A$ is unital if and only if $X$ is compact.
\end{theorem}

This actually defines a contravariant functor from the category of compact Hausdorff spaces with continuous maps to the category of unital $\mathrm{C}^*$-algebras and unital $^*$-homomorphisms.
This is the reason why $\mathrm{C}^*$-algebraist like to say that $\mathrm{C}^*$-algebras are non-commutative topological spaces. In fact, this analogy can be pushed a bit further as many topological properties are equivalent to $\mathrm{C}^*$-algebraic properties of the corresponding commutative algebra: unital - compact, non-unital - locally compact but not
compact, closed ideal - open subset, unitisation - compactification, etc. 

\section{Nuclear dimension}

This section is devoted to explaining an approach for defining a regularity condition that was deeply investigated by Winter, Kirchberg and Zacharias  \cite{Wi03,cov-dim2,Kir-Win,WZ10}. Remarkably, this notion became a driving force in the classification programme of simple separable $\mathrm{C}^*$-algebras.

The idea behind this notion is simple: since we can view $\mathrm{C}^*$-algebras as non-commutative topological spaces, we can define a notion of dimension for $\mathrm{C}^*$-algebras using topological dimension as a model. The 
first step at implementing this approach is identifying the notion of topological dimension that is most convenient for our purposes. Let us focus on this first.

\subsection{Topological dimension}

Defining the dimension of topological spaces is a rather old problem and it goes back to work of Brouwer, Lebesgue and Poincar\'e. There exist different notions of dimension, for instance \textit{small inductive dimension} (introduced independently by Urysohn and Menger \cite{Ury,Menger}) and \textit{large inductive dimension} (introduced by \v{C}ech \cite{Cech31}). We will focus mostly in another notion introduced by \v{C}ech in \cite{Cech33}, \textit{covering dimension}, that was motivated by previous work of Lebesgue \cite{Lebesgue}. This notion is more suitable for our purposes but it turns out that these three different theories agree on metrisable spaces.
We refer to \cite{Pe75} and \cite{Eng} for more comprehensive reviews of these notions.

Consider a finite open cover $\mathcal{U}$ of a topological space $X$ and a natural number $n$. The \textit{order} of $\mathcal{U}$ is at most $n$ if any point of $X$ is contained in at most $n+1$ open sets of $\mathcal{U}$. 

\begin{definition}
	Let $X$ be a topological space. The \textit{covering dimension} of $X$ is at most $n$, denoted by $\dim X \leq n$, if any finite open cover of $X$ has an open refinement of order at most $n$. 
	The covering dimension of $X$ is the minimum $n$ for which $\dim X \leq n$.

\end{definition} 

\begin{figure}[ht!]
	\centering
	\includegraphics[width=88mm]{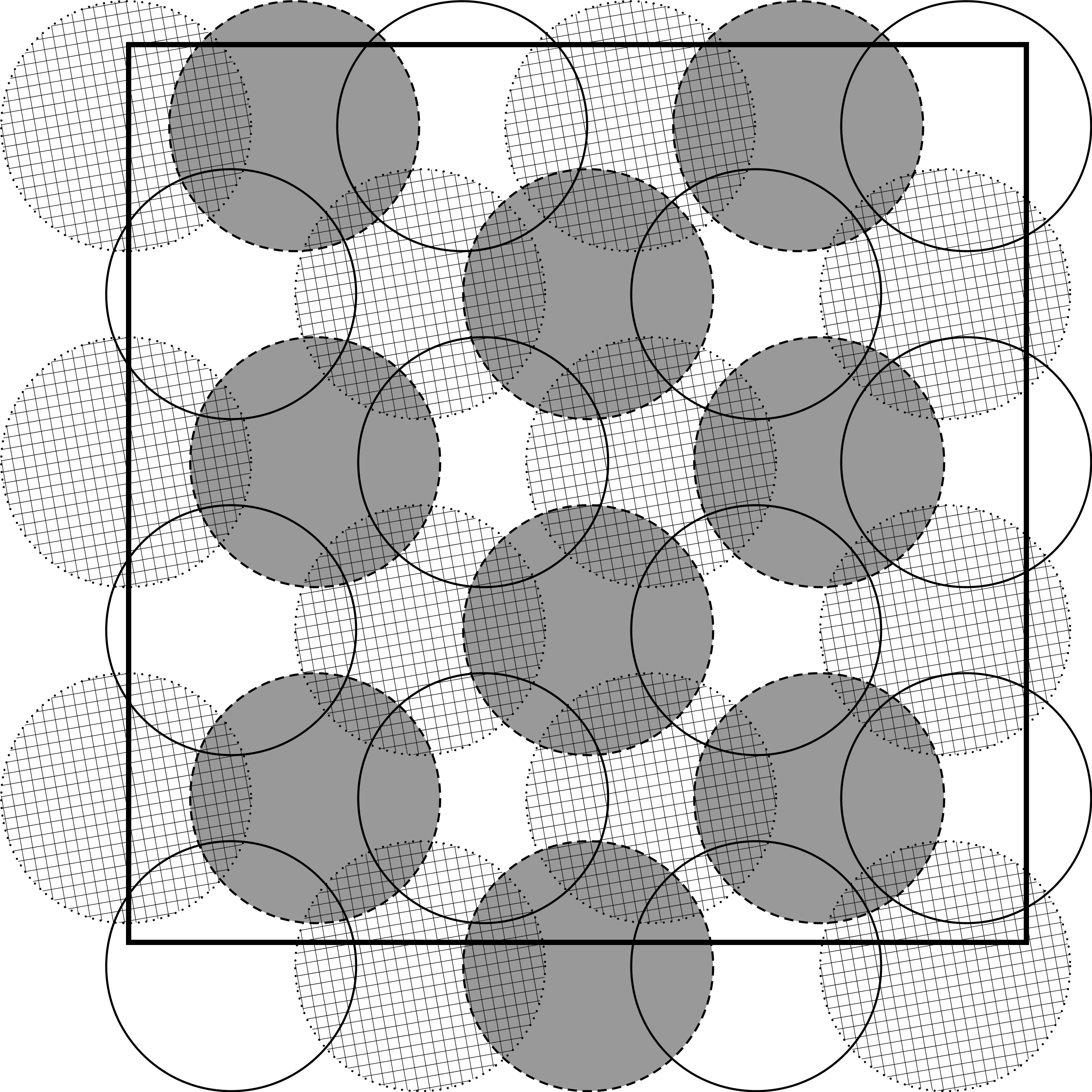}
	\caption{Covering dimension of the square is at most $2$ \label{dim.square}}
\end{figure}

This notion agrees with our intuition that the topological dimension of $\mathbb{R}^n$ must be equal to $n$ \cite[Theorem 3.2.7]{Pe75}, even though the proof of this seemingly obvious fact is not easy! (see Figure \ref{dim.square}).
Another interesting fact about covering dimension is the following:
A compact Hausdorff space has covering dimension equal to zero if and only if it is \textit{totally disconnected} \cite[Proposition 3.1.3]{Pe75}. For instance, the Cantor set has covering dimension equal to zero.

When a topological space $X$ is normal, a stronger form of covering dimension can be obtained: $\dim X \leq n$ if and only if for any finite open cover $\mathcal{V}$ there exists an open finite refinement $\mathcal{U}$ that can be decomposed into $n+1$ pairwise disjoint families $\mathcal{U}_0, \ldots, \mathcal{U}_n$ such that if $U, V$ are elements of the same $\mathcal{U}_i$ then $U \cap V = \emptyset$ \cite[Proposition 1.5]{Kir-Win}. 
This characterisation provides a nice interpretation of covering dimension in terms of colourings,
\textit{i.e.}\ we can colour the refinement $\mathcal{U}$ using $n+1$ colours, assigning the colour $i$ to each element of $\mathcal{U}_i$. The key fact of this colouring is that any two open sets with the same colour do not intersect each other.

This is the exact form of topological dimension that will be used as a model for nuclear dimension. Roughly speaking, we will identify the coloured finite open covers produced by covering dimension with certain kind of finite dimensional approximations.

\subsection{Completely positive maps of order zero}
  
Let us introduce all the ingredients we need for constructing the finite dimensional approximations induced by covering dimension.

An element in a $\mathrm{C}^*$-algebra is \textit{positive} if it is self-adjoint and its spectrum is contained in the positive real numbers. The set of positive elements of $\A$ is denoted by $\A_+$.
A linear map $\varphi: \A \to \B$ is \textit{positive} if it maps positive elements into positive elements, \textit{i.e.}\ $\varphi(\A_+) \subseteq \B_+$. 
This is a much broader class of maps between $\mathrm{C}^*$-algebras than $^*$-homomorphisms. For instance, there are no non-zero $^*$-homomorphisms between $M_2(\mathbb{C})$ and $\mathbb{C}$ but there are infinitely many positive maps from $M_2(\mathbb{C})$ to $\mathbb{C}$. An easy example can be the \textit{unnormalised trace} $\tr: M_2(\mathbb{C}) \to \mathbb{C}$ given by 
\begin{equation}
\tr \left( \left[ \begin{matrix}
a_{11} & a_{12} \\
a_{21} & a_{22}
\end{matrix} \right] \right) = a_{11} + a_{22}.
\end{equation}

We will denote by $M_n(\A)$ to the $\mathrm{C}^*$-algebra of $n \times n$ matrices over $\A$. This algebra can also be seen as the tensor product $M_n (\mathbb{C})\otimes \A$.

\begin{definition}
	A linear map $\varphi: \A \to \B$ is \textit{completely positive} if, for all $n \in \mathbb{N}$, the map $\varphi^{(n)}: M_n(\A) \to M_n(\B)$, given by
	\begin{equation}
	\varphi^{(n)} ([a_{ij}]_{i,j}) = [\varphi(a_{ij})]_{i,j},
	\end{equation}
	if positive.
	For brevity, we will use c.p.\ for  completely positive and c.p.c.\ for completely positive and contractive.
\end{definition}

If $\varphi$ is a $^*$-homomorphism, then each $\varphi^{(n)}$ is a $^*$-homomorphism as well. Hence the $^*$-homomorphisms $\varphi$ is a c.p.\ map. However, the converse is not true. For instance, let $p \in \A$ be a non trivial \textit{projection} (i.e.\  $p=p^*=p^2$), then it is straightforward to verify that the map $a \mapsto pap$ is c.p.c.\ but it is not a $^*$-homomorphism in general.

Thanks to a theorem due to Stinespring \cite[Theorem 1]{Stin}, we know the structure of c.p.\ maps: A bounded linear map $\varphi: \A \to B(\H)$ is c.p.\ if and only if there exists a Hilbert space $\widetilde{\H}$, a $^*$-homomorphism $\pi: \A \to B(\widetilde{\H})$, and an operator $V: H \to \widetilde{\H}$ such that
\begin{equation}
\varphi(a) = V^* \pi(a) V, \qquad a \in \A.
\end{equation}

Completely positive maps are an essential ingredient in the approximation theory of $\mathrm{C}^*$-algebras and we refer to \cite{Pau02} for a detailed introduction to c.p.\ maps. However, this class of maps does not preserve the structure enough for some of our purposes.
The right class of maps was introduced by Winter when the domain is finite dimensional and by Winter and Zacharias in general \cite{Wi03, WZ09}.

\begin{definition}
	A c.p.\ map $\varphi: \A \to \B$ is \textit{order zero} if it preserves orthogonality, \textit{i.e}.\ $\varphi(a)\varphi(b)=0$ if $a,b \in \A$ satisfy $ab=0$.  
\end{definition}

This class of maps lies between $^*$-homomorphisms and c.p.\ maps. Indeed, c.p.\ maps of the form $a \mapsto pap$ are not order zero in general and  order zero maps of the form $a \mapsto \lambda a$, with $\lambda \in (0,1)$, are not $^*$-homomorphisms.
As before, we have a theorem \'{a} la Stinespring that unravels the structure of order zero maps.

\begin{theorem}[\cite{WZ09}]\label{thm:structure.order.zero}
	Let $\varphi: \A \to \B$ be an order zero map. Then there exists a unital $\mathrm{C}^*$-algebra $\mathcal{M}$ containing $\varphi(\A)$,\footnote{A precise description of the algebra $\mathcal{M}$ can be found in \cite{WZ09}.} 
	$h \in \mathcal{M}_+$ commuting with $\varphi(A)$, and a $^*$-homomorphism $\pi: \A \to \mathcal{M}$ such that
	\begin{equation}
	\varphi(a) = h\pi(a) = \pi(a) h, \qquad a \in \A.
	\end{equation}
	If $\A$ is unital, then $h = \varphi(1_{\A})$.
\end{theorem}

This theorem enables us to obtain a canonical bijection between c.p.c.\ order zero maps $\A \to \B$ and $^*$-homomorphisms $C_0((0,1], \A) \to \B$ \cite[Corollary 4.1]{WZ09}.
The $\mathrm{C}^*$-algebra $C_0((0,1], \A)$ is called the \textit{cone over $\A$} and it is given by
\begin{equation}
C_0((0,1], \A) := \{ g: (0,1] \to \A \mid g \text{ is continuous and} \lim_{t \to 0} g(t) =0 \}.
\end{equation}
This algebra is generated by the functions $f\otimes a: (0,1] \to \A$ given by 
\begin{equation}
(f\otimes a) (t)= f(t)a,
\end{equation}
where $f \in C_0(0,1], \ a \in A$ and $t \in (0,1]$.

The bijection is precisely given in the following way: if $\varphi: \A \to \B$ is an order zero map, then the map $\rho_\varphi: C_0((0,1], \A) \to \B$ given on elementary tensors by 
\begin{equation}
\rho_\varphi(f\otimes a) := f(h)\pi(a)
\end{equation}
is a $^*$-homomorphism. Conversely, if $\rho: C_0((0,1], \A) \to \B$ is a $^*$-homomorphism, then the map $\varphi_\rho: \A \to \B$ given by
\begin{equation}
\varphi_\rho (a) := \rho(\mathrm{id}_{(0,1]} \otimes a)
\end{equation} 
is a c.p.c.\ order zero map. 
This is one of the key features of order zero maps and also explains why order zero maps are so useful. They avoid \textit{topological obstructions} in contrast to the situation with $^*$-homomorphisms, where due to topological obstructions, non-trivial $^*$-homomorphisms might not exist.
Indeed, the $\mathrm{C}^*$-algebra $C_0((0,1], \A)$ is homotopy equivalent to the zero  $\mathrm{C}^*$-algebra.
Hence, order zero maps do not carry any topological information.

\subsection{Covering dimension of $C(X)$}\label{subsection:dim.C(X)}

Before defining nuclear dimension, we will discuss the commutative case as the motivating example for the proper definition. 
In order to simplify our equations, we will write $a \approx_{\epsilon} b$ instead of $\|a - b \|< \epsilon$.

Consider a second countable compact Hausdorff space $X$ of covering dimension equal to $n$, $\epsilon > 0$ and let $\mathfrak{F}$ be a finite set of $C(X)$. 
By compactness of $X$, we can find a finite open cover $\mathcal{V}$ of $X$ such that if $x,y \in V$, for some $V \in \mathcal{V}$, then $f(x) \approx_{\epsilon} f(y)$ for all $f \in \mathfrak{F}$.

Using the covering dimension of $X$, we can find a finite open refinement $\mathcal{U}$ of order $n+1$. So we can decompose it in the following way: $\mathcal{U} = \mathcal{U}_0 \sqcup \ldots \sqcup \mathcal{U}_n$ and $U \cap V = \emptyset$ if $U$ and $V$ are elements $\mathcal{U}_k$ for some $k = 0, \ldots, n$.
Fix points $x_U \in U$ for each $U \in \mathcal{U}$ and set a c.p.c.\ map $\psi: C(X) \to \mathbb{C}^{\mathcal{U}}$ by 
\begin{equation}
\psi(f) := (f(x_U))_{U \in \mathcal{U}}, \qquad f \in C(X).
\end{equation} 

Let $(h_U)_{U \in \mathcal{U}}$ be a \textit{partition of unity} subordinated to the cover $\mathcal{U}$. Recall that this means that $(h_U)$ is a family of continuous functions $X \to \mathbb{R}_+$ such that the support of $h_U$ is contained in $U$ and  $\sum_{U \in \mathcal{U}}h_U (x) =1$ for all $x \in X$. Set a c.p.\ map $\varphi: \mathbb{C}^{\mathcal{U}} \to C(X)$ by 
\begin{equation}
\varphi((x_U)) := \sum_{U \in \mathcal{U}} x_U h_U, \qquad (x_U) \in \mathbb{C}^{\mathcal{U}}.
\end{equation} 
It is straightforward to obtain
\begin{equation}
\varphi \circ \psi (f)(y) = \sum_{U \in \mathcal{U}} f(x_U) h_U (y) = \sum_{y \in U} f(x_U) h_U (y) \approx_{\epsilon} f(y).
\end{equation} 

This is the standard proof of the fact that $C(X)$ is a \textit{nuclear} $\mathrm{C}^*$-algebra (see \cite[Proposition 2.4.2]{BO08}). However, this approximation we have constructed has more properties. Let $\varphi_k: \mathbb{C}^{\mathcal{U}} \to C(X)$ be the c.p.\ map given by restricting $\varphi$ to the family $\mathcal{U}_k$, i.e.\
\begin{equation}
\varphi_k ((x_U)_{U \in \mathcal{U}}) := \sum_{U \in \mathcal{U}_k} x_U h_U.
\end{equation} 
If we view $\mathbb{C}^{\mathcal{U}} = \bigoplus_{k=0}^n \mathbb{C}^{\mathcal{U}_k}$, then $\varphi_k = \varphi|_{\mathbb{C}^{\mathcal{U}_k}}$. Each $\varphi_k$ is a c.p.c.\ order zero map and this follows from the fact that any two open sets in $\mathcal{U}_k$ are disjoint. This is precisely where the extra topological properties of the refinement $\mathcal{U}$ come into play. Therefore, the c.p.\ map $\varphi$ is the sum of $n+1$ c.p.c.\ order zero maps,
\begin{equation}
\varphi = \sum_{k=0}^n \varphi_k.
\end{equation}

The nuclear dimension of $C(X)$ will be defined as the number of order zero summands that the second map $\varphi$ has minus one. The preceding discussion shows that the non-commutative covering dimension of $C(X)$ is at most $n$.
The definition of nuclear dimension mimics the approximations we have described above for $C(X)$ and this is exactly why we refer to it as a non-commutative covering dimension.

\subsection{The nuclear dimension of a $\mathrm{C}^*$-algebra}\label{subsection:dim.nuc}
After discussing the commutative case, let us summarise the key properties of the completely positive approximations constructed in this case:
\begin{enumerate}
	\item  we can approximate the identity map $\mathrm{id}_{C(X)}$ with a finite dimensional algebra $\mathbb{C}^r$ and maps $ C(X) \xrightarrow{\psi} \mathbb{C}^r \xrightarrow{\varphi} C(X)$, i.e.\ for a finite subset of $C(X)$ and up to some small error, $\mathrm{id}_{C(X)} \approx \varphi \psi$,
	\item the map $\psi: C(X) \to \mathbb{C}^r$ is c.p.c., 
	\item there is a decomposition of $\mathbb{C}^r$ into $n+1$ finite dimensional subalgebras, $\mathbb{C}^r = \bigoplus_{k=0}^n \mathbb{C}^{m_i}$, so that each restriction $\varphi|_{\mathbb{C}^{m_i}}$ is a c.p.c.\ order zero map. Thus $\varphi$ can be expressed as the sum of $n+1$ order zero maps.
\end{enumerate}

This leads us to the actual definition of nuclear dimension introduced by Winter and Zacharias, where they replaced the commutative algebras with more general ones.

\begin{definition}[\cite{WZ10}]\label{def:dim.nuc}
	A $\mathrm{C}^*$-algebra $\A$ has \textit{nuclear dimension} at most $n$, $\dim_{\mathrm{nuc}} A \leq n$, if for any finite set $\mathfrak{F} \subseteq A$ and $\epsilon > 0$, there exist finite dimensional $\mathrm{C}^*$-algebras $F_0, F_1 \ldots, F_n$ and maps $\psi:\A \to \bigoplus_{k=0}^n F_k$, $\varphi: \bigoplus_{k=0}^n F_k \to \A$ such that
	\begin{enumerate}[(i)]
		\item $\|a - \varphi \circ \psi (a) \| < \epsilon$ for all $a \in \mathfrak{F}$,
		\item $\psi$ is c.p.c.,
		\item $\varphi_k := \varphi|_{F_k}$ is a c.p.c.\ map of order zero for $i=0, \ldots, n$.
	\end{enumerate} 
\end{definition}

The name `nuclear dimension' comes from the fact that $\mathrm{C}^*$-algebras with approximations of the form $\A \xrightarrow{\psi} F \xrightarrow{\varphi} \A$, with $F$ finite dimensional and $\psi, \varphi$ are c.p.c.\ maps, are called \textit{nuclear}, see \cite[Chapter 2]{BO08} for more information about nuclearity. It follows that $\mathrm{C}^*$-algebras with finite nuclear dimension are nuclear but the converse is false, for instance $C([0,1]^{\mathbb{N}})$.

An interesting variant of the previous definition is \textit{decomposition rank}, denoted by $\mathrm{dr}\, (\A)$. This precursor notion, introduced by Winter and Kirchberg in \cite{Kir-Win}, is very similar to nuclear dimension but it asks the second map $\varphi$ to be contractive.
Observe that in the definition of nuclear dimension we only know that the norm of $\varphi$ is at most $n+1$. 
This seemingly insubstantial change is rather subtle and deep since it imposes stronger structural properties (see Section \ref{section:Cuntz.alg}). 
Note that we always have 
\begin{equation}
\dim_{\mathrm{nuc}}\, \A \leq \mathrm{dr}\, \A.
\end{equation}
We refer to \cite{WiSurvey} for a more detailed explanation. 

After the discussion in the previous subsection, the following theorem should not be a great surprise. 

\begin{theorem}[{\cite[Proposition 2.4]{WZ10}}]\label{thm:dim.C(X)}
	Let $X$ be a compact Hausdorff space. Then
	\begin{equation}
	\dim_{\mathrm{nuc}} C(X) = \dim X.
	\end{equation}
\end{theorem}

We already showed that $\dim_{\mathrm{nuc}} C(X) \leq \dim X$ but the other direction is a bit more subtle. In Section \ref{ex.commutative}, we will present a direct proof of this fact.

We finish this section by pointing out that the notion of nuclear dimension is a well-behaved dimension theory and it has good permanence properties. For instance,
\begin{align}
\dim_{\mathrm{nuc}} \A \oplus \B =& \max \{ \dim_{\mathrm{nuc}} \A, \dim_{\mathrm{nuc}} \B\}, \notag \\
\dim_{\mathrm{nuc}} \A \otimes \B \leq & \left( \dim_{\mathrm{nuc}} \A + 1\right)\left( \dim_{\mathrm{nuc}} \B + 1\right)-1, \notag \\
\max \{\dim_{\mathrm{nuc}} I, \dim_{\mathrm{nuc}} \A / I \} \leq& \dim_{\mathrm{nuc}} \A \leq \dim_{\mathrm{nuc}} I + \dim_{\mathrm{nuc}} \A / I + 1,
\end{align}
where in the last inequality, $I$ is a closed two sided ideal of $\A$ \cite[Proposition 2.3]{WZ10}.

In the following sections of this survey, we present more examples of $\mathrm{C}^*$-algebras and discuss recent developments regarding their nuclear dimension.

\section{Zero dimensional objects}
	A natural question regarding nuclear dimension is to identify the $\mathrm{C}^*$-algebras with nuclear dimension exactly equal to zero. We already mentioned that a compact Hausdorff space has covering dimension zero if and only if it is totally disconnected. In conjunction with Theorem \ref{thm:dim.C(X)}, this immediately provides us a characterisation of unital commutative $\mathrm{C}^*$-algebras with nuclear dimension equal to zero: $C(X)$ has nuclear dimension zero if and only if $X$ is totally disconnected. 
	
	In the non-commutative setting, it is also straightforward to see that finite dimensional $\mathrm{C}^*$-algebras have nuclear dimension equal to zero, but we would also like to know if there exist infinite dimensional non-commutative $\mathrm{C}^*$-algebras with nuclear dimension equal to zero. 
	The answer is precisely those $\mathrm{C}^*$-algebras that locally ``look like'' finite dimensional $\mathrm{C}^*$-algebras.
	
	\begin{definition}
	A separable $\mathrm{C}^*$-algebra $\A$ is \textit{approximately finite dimensional} (AF) if there exists an increasing sequence of finite dimensional $\mathrm{C}^*$-subalgebras $(F_n)_{n=1}^\infty$ such that 
	\begin{equation}
	\A = \overline{\bigcup_{n=1}^\infty F_n}. 
	\end{equation}	
	\end{definition}

	Bratteli showed in \cite[Theorem 2.2]{Bra} that a separable $\mathrm{C}^*$-algebra $\A$ is AF if and only if for any finite subset $\mathfrak{F} \subseteq \A$ and $\epsilon>0$ there exists a finite dimensional $\mathrm{C}^*$-subalgebra $F \subseteq \A$ such that
	\begin{equation}
	\mathrm{dist}(a, F) : = \inf_{x \in F} \| a - x  \| < \epsilon, \qquad a \in \mathfrak{F}.
	\end{equation}
	This is a local characterisation of AF-algebras that is very useful when computing their nuclear dimension.
	
	Before going further, let us consider some examples of an AF-algebra.
	The $\mathrm{C}^*$-algebra of compact operators on a separable Hilbert space $\mathbb{K}(\H)$ is a (non-unital) AF-algebra. Indeed, suppose $(e_n)_{n=1}^\infty$ is an orthonormal basis of $\H$. The projection $e_i \otimes e_i : \H \to \mathbb{C} e_i \subseteq \H$ is given by 
	\begin{equation}
	e_i \otimes e_i (h) := \langle h, e_i \rangle e_i, \qquad h \in \H.
	\end{equation}
	Let $p_n$ be the projection $\sum_{i=1}^n e_i \otimes e_i$. It is not difficult to see that 
	\begin{equation}
	\mathbb{K}(\H) = \overline{\bigcup_{n=1}^\infty p_n \mathbb{K}(\H) p_n} \quad \text{and} \quad  p_n \mathbb{K}(\H) p_n \cong M_n(\mathbb{C}).
	\end{equation}
	This shows that $\mathbb{K}(\H)$ is an AF-algebra.
	
	Another example of an AF-algebra is the algebra of continuous functions on the Cantor set. In fact, 
	any commutative $\mathrm{C}^*$-algebra $C_0(X)$ is AF if and only if $X$ is totally disconnected; this shows that the commutative $\mathrm{C}^*$-algebras with nuclear dimension equal to zero that we already identified are indeed AF-algebras.
	
	Winter showed that the zero dimensional objects for nuclear dimension are precisely the AF-algebras. This essentially relies on the local characterisation of AF-algebras that we discussed above.
	
	\begin{theorem}[\cite{WZ10}]
		A separable $\mathrm{C}^*$-algebra $\A$ is approximately finite dimensional if and only if $\dim_{\mathrm{nuc}} \A = 0$.
	\end{theorem}

\section{$\mathrm{C}^*$-algebras associated to groups}

In this section, we will explain how to construct $\mathrm{C}^*$-algebras associated to discrete groups.
Let us recall some basic definitions first.

An operator $U \in B(\H)$ is a \textit{unitary} if $U^* U = U U^* = 1_{\H}$. We denote by $\mathcal{U}(\H)$ to the set of unitaries of $B(\H)$. This set has the structure of a group. Indeed, 
by definition, the inverse of a unitary is precisely its adjoint, which is a unitary as well. Note that $\mathcal{U}(\H)$ is also closed under multiplication, so it is a group.

The \textit{group ring} $\mathbb{C} \Gamma$ is given by
\begin{equation}
\mathbb{C} \Gamma = \biggl\{ \sum_{g \in \Gamma} a_g \cdot g \, \mathrel{\Big|} \, a_g \in \mathbb{C} \text{ and } a_g =0 \text{ except for finitely many} \biggl\}.
\end{equation}
A \textit{unitary representation} of a discrete group $\Gamma$ on a Hilbert space $\H$ is a group homomorphism $\pi: \Gamma \to \mathcal{U}(\H)$. It follows that a unitary representation of $\Gamma$ extends by linearity to a representation of $\mathbb{C} \Gamma$ into $B(\H)$.
Naturally, a unitary representation of $\Gamma$ induces a seminorm on $\mathbb{C} \Gamma$ by simply taking 
\begin{equation}
\| x \|_\pi : = \| \pi(x) \|_{B(\H)}, \qquad x \in \mathbb{C} \Gamma.
\end{equation}

An important example of a unitary representation of a discrete group $\Gamma$ is the \textit{left regular representation}. Consider the Hilbert space 
\begin{equation}
\ell^2 (\Gamma): = \biggl\{ a: \Gamma \to \mathbb{C} \ \mathrel{\Big|} \ \sum_{g \in \Gamma} | a_g |^2 < \infty \biggl\},
\end{equation}
and define $\lambda: \Gamma \to B(\ell^2 (\Gamma))$ by
\begin{equation}
\lambda_g (a)(h) = a(g^{-1} h), \qquad g,h \in \Gamma, \ a \in \ell^2(\Gamma)
\end{equation}
The seminorm $\|\cdot\|_\lambda$ is in fact a norm since the vector $\delta_e \in \ell^2$ separates points of $\lambda(\mathbb{C} \Gamma)$.

The \textit{reduced group $\mathrm{C}^*$-algebra of $\Gamma$}, denoted by $\mathrm{C}^*_r(\Gamma)$, is given by completing $\mathbb{C} \Gamma$ with respect to the $\mathrm{C}^*$-norm $\|\cdot \|_\lambda$ induced by the left regular representation $\lambda$,  
\begin{equation}
\mathrm{C}^*_r(\Gamma) := \overline{ \mathbb{C} \Gamma}^{\|\cdot \|_\lambda}.
\end{equation}

Another norm that we can consider is the \textit{universal norm} given by
\begin{equation}
\| x \|_u = \sup_{\pi} \| x \|_\pi, \qquad x \in \mathbb{C} \Gamma,
\end{equation} 
where the supremum is taken over all unitary representations $\pi$ of $\Gamma$. The \textit{universal group $\mathrm{C}^*$-algebra of $\Gamma$}, denoted by $\mathrm{C}^*(\Gamma)$, is given by completing $\mathbb{C} \Gamma$ with respect to the universal norm, i.e.\
\begin{equation}
\mathrm{C}^*(\Gamma) := \overline{ \mathbb{C} \Gamma}^{\|\cdot \|_u}.
\end{equation}

In particular, if the group $\Gamma$ is finite then the reduced and universal group $\mathrm{C}^*$-algebras are finite dimensional, and in fact, they are equal to the group ring. This can be deduced from the fact that $\mathbb{C} \Gamma$ is finite dimensional and the uniqueness of the $\mathrm{C}^*$-norm.

Another interesting example is when $\Gamma$ is commutative. In this case, again the reduced and universal group $\mathrm{C}^*$-algebras agree. But more is true: they are isomorphic to the $\mathrm{C}^*$-algebra of continuous functions on its \textit{Pontryagin dual} $\widehat{\Gamma}$,
\begin{equation}
\mathrm{C}^*(\Gamma) \cong \mathrm{C}^*_r(\Gamma) \cong C(\widehat{\Gamma}).
\end{equation}
This immediately tells us that, for instance, $\mathrm{C}^*(\mathbb{Z}) \cong \mathrm{C}^*_r(\mathbb{Z}) \cong C(\mathbb{T})$.

The previous examples show that, in some cases, reduced and universal group $\mathrm{C}^*$-algebras agree.
An interesting and natural question is for which class these two group $\mathrm{C}^*$-algebras agree. 
The answer is related to the famous \textit{Banach-Tarsky paradox} that states that the sphere $\mathbb{S}^2 \subseteq \mathbb{R}^3$ can be decomposed into finitely many pieces that can be reassembled as two copies of $\mathbb{S}^2$. 
A similar statement is true for spheres in higher dimensions but interestingly this cannot be done for $\mathbb{S}^1 \subseteq \mathbb{R}^2$.  
Despite the apparently geometric nature of the paradox, this is truly a measure theoretic phenomenon. This was identified by von Neumann, who introduced the concept of amenability \cite{vNe29}. 

\begin{definition}
	A discrete group $\Gamma$ is \textit{amenable} if there exists a a \textit{finitely
	additive measure} $\mu$ on $\Gamma$ that satisfies $\mu(\Gamma) = 1$ and $\mu(gA) = \mu(A)$ for $g\in \Gamma$ and $A \subseteq \Gamma$.
\end{definition}

The reason behind the Banach-Tarsky paradox for $\mathbb{S}^2$ lies in the fact that the group of isometries of $\mathbb{R}^3$ is not amenable. Meanwhile, this paradox does not hold for $\mathbb{S}^1$ precisely because the group of isometries of $\mathbb{R}^2$ is indeed amenable.

Examples of amenable discrete groups can be found in the finite groups, abelian groups, virtually solvable groups and finitely generated groups of sub-exponential growth. On the other hand, the free groups $\mathbb{F}_n$ are the canonical examples of non amenable groups. Moreover, any group that contains a copy of a free group is not amenable.

There exist a vast amount of characterisations of amenable groups (see \cite{Pat88}). Let us present two of them that relate amenability with properties of group $\mathrm{C}^*$-algebras. 

\begin{theorem}[{\cite[Theorem 2.6.8]{BO08}}]
	Let $\Gamma$ be a discrete group. The following are equivalent:
	\begin{enumerate}[(i)]
		\item $\Gamma$ is amenable,
		\item $\mathrm{C}^*_r(\Gamma) \cong \mathrm{C}^*(\Gamma)$,
		\item $\mathrm{C}_r^*(\Gamma)$ is a nuclear $\mathrm{C}^*$-algebra.
	\end{enumerate}
\end{theorem}

This theorem immediately shows that if $\Gamma$ is not amenable, its reduced group $\mathrm{C}^*$-algebra $\mathrm{C}^*_r(\Gamma)$ is not a nuclear $\mathrm{C}^*$-algebra and hence $\dim_{\mathrm{nuc}} \mathrm{C}_r^*(\Gamma) = \infty$.
	
What can we say about the nuclear dimension of $\mathrm{C}^*$-algebras associated to amenable groups? First of all, we already know that if $\Gamma$ is commutative, its group $\mathrm{C}^*$-algebra is isomorphic to $C(\widehat{\Gamma})$. By Theorem \ref{thm:dim.C(X)}, $\dim_{\mathrm{nuc}} \mathrm{C}_r^*(\Gamma) = \dim \widehat{\Gamma}$.
As a consequence, we obtain examples of amenable groups with infinite nuclear dimension, for instance $\mathbb{Z}^{\mathbb{Z}}$.
Finding amenable groups with infinite nuclear dimension appears harder
if we restrict to finitely generated groups. 
However, such examples do exist: the nuclear dimension of $\mathrm{C}^*_r(\mathbb{Z} \wr \mathbb{Z})$ is infinite \cite{GK10, Win12}.
We finish our discussion of groups and their nuclear dimension by presenting a class of amenable groups with finite nuclear dimension.

\begin{theorem}[{\cite{EGM}}]
	Let $\Gamma$ be a finitely generated \textit{virtually nilpotent}  group. Then $\dim_{\mathrm{nuc}} \mathrm{C}^*_r(\Gamma) < \infty$.
\end{theorem}

In the previous theorem, the upper bound depends on the \textit{Hirsch number} of $\Gamma$.

\section{$\mathrm{C}^*$-algebras associated to dynamical systems}
	
	Let $X$ be a compact Hausdorff space and consider a homeomorphism $\alpha: X \to X$. This homeomorphism induces an action of $\mathbb{Z}$ into $C(X)$, that we will also denote by $\alpha$, in the following way:
	\begin{equation}
	\alpha_n (f) = f \circ \alpha^{-n}, \qquad f \in C(X), \, n \in \mathbb{Z}.
	\end{equation}
	This action gives rise to a $\mathrm{C}^*$-algebra $C(X) \rtimes_{\alpha} \mathbb{Z}$ contained in $B(L^2(X) \otimes \ell^2(\mathbb{Z})) $ that is generated by a copy of $C(X)$ and a unitary $u \in B(L^2(X) \otimes \ell^2(\mathbb{Z}))$ that encodes the action, i.e.\
	\begin{equation}
	u f u^* = f \circ \alpha^{-1}, \qquad f \in C(X).
	\end{equation}
	
	One class of interesting examples of this construction are the \textit{rotation algebras} $\A_\theta$. Let $\mathbb{T} \subseteq \mathbb{C}$ be the unit circle in the complex plane. For $\theta \in [0,1)$, define $r_\theta: \mathbb{T} \to \mathbb{T}$ by $r_\theta(z) := e^{2\pi i\theta} z$.
	The corresponding crossed product is denoted by $\A_\theta := C(\mathbb{T}) \rtimes_{r_\theta} \mathbb{Z}$. 
	These algebras are called rotation algebras and are also commonly referred as non-commutative tori. This name is motivated by the fact that $\A_0 \cong C(\mathbb{T}^2)$. 
	
	These algebras have received much attention since they are some of the most natural examples of non-commutative manifolds \cite{Rie90}. For instance,
	$\A_\theta$ is isomorphic to $\A_\eta$ if and only if $\theta = \eta$ or $\theta = 1- \eta$ \cite{PV80, Rie81,Yi86}. Furthermore, if $\theta$ is rational, $\A_\theta$ is a subalgebra of $M_k (C(\mathbb{T}^2))$ for some $k \in \mathbb{N}$ \cite{Rie83}. On the other hand, if $\theta$ is irrational, $\A_\theta$ is a simple $\mathrm{C}^*$-algebra and locally ``looks like" $C(\mathbb{T})\otimes F$ where $F$ is a finite dimensional $\mathrm{C}^*$-algebra \cite{EV93}. 
	
	This shows that the nuclear dimension of the rotation algebra $\A_\theta$ is one when $\theta$ is irrational, or two when $\theta$ is rational. Indeed, this essentially follows from the fact that $\dim_{\mathrm{nuc}} M_k (\mathbb{T}^2) = \dim \mathbb{T}^2 = 2$, and $\dim_{\mathrm{nuc}} C(\mathbb{T})\otimes F = \dim \mathbb{T} = 1$ \cite[Example 6.1]{WZ10}.
	
	In general, understanding the internal structure of $C(X) \rtimes_{\alpha} \mathbb{Z}$ can be a very challenging task. However,  
	it was recently shown in \cite{HW17} that when the covering dimension of $X$ is finite, the nuclear dimension of the $\mathrm{C}^*$-algebra induced by any homeomorphism is finite as well. Precisely
	\begin{equation}\label{eq:nd.gen.homeo}
	\dim_{\mathrm{nuc}} C (X) \rtimes_{\alpha} \mathbb{Z} \leq 2 (\dim X)^2 + 6 \dim X + 4.
	\end{equation}
	This is a surprising result because it does not impose conditions on the homeomorphism or on the space $X$ (apart from having finite covering dimension).

	The $\mathrm{C}^*$-algebra $C(X) \rtimes_{\alpha} \mathbb{Z}$ is an example of a more general construction. Given a $\mathrm{C}^*$-algebra $\A \subseteq B(\H)$,  a discrete group $\Gamma$ and a group action $\alpha: \Gamma \curvearrowright A$, we can always form the \textit{reduced crossed product} $\mathrm{C}^*$-algebra $\A \rtimes_{\alpha,r} \Gamma$ that is contained in $B(\H \otimes \ell^2(\Gamma))$. In general, determining the structure of crossed products can be a very difficult task but there exist recent spectacular results that enable us to estimate the nuclear dimension of certain crossed products \cite{Sza15,HWZ15, EN17, HSWW17, SWZ17,KS18}. 

\section{Cuntz algebras}\label{section:Cuntz.alg}

Let $\H$ be a separable infinite dimensional Hilbert space. For any natural number $n$, we can decompose $\H$ as $\H_1 \oplus ... \oplus \H_n$ with $\H_i \cong \H$ for all $i$. Then, we can find elements $v_1, \ldots, v_n \in B(\H)$ implementing an isomorphism between $\H$ and $\H_i$ via conjugation that satisfy 
\begin{align}
\sum_{i=1}^{n} v_i v_i^* = 1_{B(\H)} \quad \text{and} \quad v_1^* v_1 = v_2^* v_2 = \ldots = v_n^* v_n = 1_{B(\H)}.
\end{align}
The \textit{Cuntz algebra} $\mathcal{O}_n$ is the $\mathrm{C}^*$-subalgebra of $B(\H)$ generated by such elements $v_1, \ldots, v_n$. 
Similarly, we can find a sequence $(u_i)_{i \in \mathbb{N}}$ of elements in $B(\H)$ satisfying
\begin{equation}
u_i^* u_i = 1_{B(\H)}, \qquad (u_i u_i^*) (u_j u_j^*) = 0, \qquad i\neq j, \quad  i,j \in \mathbb{N}.
\end{equation}
The \textit{Cuntz algebra} $\mathcal{O}_\infty$ is the $\mathrm{C}^*$-subalgebra of $B(\H)$ generated by the sequence $(u_i)_{i \in \mathbb{N}}$.

These algebras were introduced by Cuntz in \cite{Cu77} and they possess very interesting properties. For instance, each $\mathcal{O}_n$ is simple, nuclear, unital, and given any two non-zero positive elements $a,b \in \mathcal{O}_n$, there exists $x \in \mathcal{O}_n$ such that $a = x^* b x$. Simple $\mathrm{C}^*$-algebras that satisfy this property are called \textit{purely infinite}. 
There exist spectacular results concerning Cuntz algebras (and purely infinite $\mathrm{C}^*$-algebras) that are beyond the scope of this survey. We refer to \cite{Ro02} for a comprehensive overview. 

We briefly mentioned the notion of \textit{decomposition rank} in Section \ref{subsection:dim.nuc}. This stronger notion requires  approximations $\A \xrightarrow{\psi} F \xrightarrow{\phi} \A$ similar to the approximations in Definition \ref{def:dim.nuc} but also requires the second map $\varphi$ to be contractive, i.e.\  $\| \varphi \| \leq 1$. 
Finite decomposition rank imposes stronger structural properties on the algebra than finite nuclear dimension, for instance \textit{strong quasidiagonality} \cite[Proposition 5.1]{Kir-Win}. Cuntz algebras are not even \textit{quasidiagonal}, and hence, they cannot have finite decomposition rank. However, these algebras do have finite nuclear dimension.

\begin{theorem}[{\cite[Theorem 7.4]{WZ10}, \cite[Theorem 4.1]{En15}}]
	For any $n \in \mathbb{N} \cup \{\infty\}$, $\dim_{\mathrm{nuc}} \mathcal{O}_n =1$. 
\end{theorem}

This theorem is a special case of a more general result. 
Cuntz algebras are particular examples of separable simple purely infinite and nuclear $\mathrm{C}^*$-algebras -- a class of $\mathrm{C}^*$-algebras commonly referred to as the \textit{Kirchberg algebras}. It was showed in \cite{BBSTWW,RSS15} that the nuclear dimension of any Kirchberg algebra is exactly 1. More generally, for any non-zero separable nuclear $\mathrm{C}^*$-algebra $\A$, the nuclear dimension of $\A \otimes \mathcal{O}_\infty$ is finite \cite{Sza17}. It was recently showed that  
$\dim_{\mathrm{nuc}} \A \otimes \mathcal{O}_\infty = 1$ \cite{BGSW}.

\section{General simple $\mathrm{C}^*$-algebras}

Now, we focus our attention on the nuclear dimension of simple separable $\mathrm{C}^*$-algebras. It turns out that for this class of $\mathrm{C}^*$-algebras, there exists an interesting trichotomy. We emphasize that commutative $\mathrm{C}^*$-algebras apart from $\mathbb{C}$ are never simple. 

\begin{theorem}[{\cite{CETWW,CE}}]
	The nuclear dimension of a simple separable $\mathrm{C}^*$-algebra is $0,1$ or $\infty$.
\end{theorem}

Simple $\mathrm{C}^*$-algebras with infinite nuclear dimension include $\mathrm{C}^*$-algebras that are not nuclear, as $\mathrm{C}^*_r(\mathbb{F}_2)$, but also the exotic Villadsen algebras \cite{Vi98,Vi99} and the famous examples constructed by R\o rdam \cite{Ror03} and Toms \cite{Toms}.
On the other hand, unital simple (UCT) $\mathrm{C}^*$-algebras with finite nuclear dimension are those that can be classified using their K-theory and traces \cite{GLN,EGLN,TWW} (see \cite{WiSurvey} for a survey on the subject). 

The reason behind this trichotomy is that finite nuclear dimension imposes strong conditions in the algebra, and one of those conditions is so strong that it actually forces the nuclear dimension to be at most one.
In order to explain this more clearly, let us introduce an exotic $\mathrm{C}^*$-algebra: the Jiang-Su algebra $\mathcal{Z}$ \cite{JS99}. This algebra can be seen as an infinite-dimensional $\mathrm{C}^*$-version of the complex numbers and it has many interesting properties. This algebra, which locally looks like some subalgebra of $C([0,1], M_{pq}(\mathbb{C}))$ where $p$ and $q$ are relatively prime numbers, is simple and separable with $\dim_{\mathrm{nuc}} \mathcal{Z} =1$. We refer to \cite{RW10} for a detailed explanation of the Jiang-Su algebra $\mathcal{Z}$.

In a notable tour de force, Winter showed that finite nuclear dimension implies tensorial absorption of the Jiang-Su algebra, $\A \otimes \mathcal{Z} \cong \A$, when the algebra is separable simple unital and non elementary\footnote{A $\mathrm{C}^*$-algebra $\A$ is elementary if it is isomorphic to $M_n(\mathbb{C})$ for some $n$, or the compacts $\mathbb{K}(\H)$.} \cite{Win12}. This result was later extended to the non-unital case by Tikuisis \cite{Tik14}. 

It was conjectured by Toms and Winter that the converse direction should hold as well and this was verified under certain extra conditions, see \cite{MS14,SWW,BBSTWW}. 
Recently, this direction was fully verified, first for
unital $\mathrm{C}^*$-algebras in \cite{CETWW}, then extended to the non-unital case
in \cite{CE}. It was actually shown that if the $\mathrm{C}^*$-algebra is $\mathcal{Z}$-stable then its nuclear dimension is at most one.

By the results mentioned above, if the nuclear dimension of a simple separable $\mathrm{C}^*$-algebra $\A$ is finite, then it has to absorb the Jiang-Su algebra $\mathcal{Z}$, and so it must have nuclear dimension at most $1$. Finally, if $\A$ is AF then its dimension is exactly zero; otherwise, it has nuclear dimension equal to one. 

The following theorem summarises all that we have said so far.

\begin{theorem}[{\cite{Win12,Tik14, CETWW, CE}}]
	Let $\A$ be a simple separable nuclear and non elementary $\mathrm{C}^*$-algebra. The following are equivalent:
	\begin{enumerate}[(i)]
		\item $\dim_{\mathrm{nuc}} \A < \infty$,
		\item $\dim_{\mathrm{nuc}} \A \leq 1$,
		\item $\A \otimes \mathcal{Z} \cong \A$.
	\end{enumerate}
\end{theorem}

\section{Nuclear dimension of commutative algebras} \label{ex.commutative}

As already explained in Section \ref{subsection:dim.C(X)}, the nuclear dimension of the $\mathrm{C}^*$-algebra $C(X)$ agrees with the covering dimension of $X$. 
In his early work
on developing non-commutative covering dimension, Winter proved that the covering dimension of a second countable locally compact Hausdorff space $X$ agrees with the \textit{completely positive rank} of $C_0(X)$ \cite{Wi03}. This notion was a precursor concept of decomposition rank and nuclear dimension. Later on, when decomposition rank and nuclear dimension were introduced, a direct proof of this fact 
was not given; instead, its proof relied on the original proof concerning completely positive rank and some relations between this former concept and nuclear dimension. 
In this section, we present a direct proof of this fact in the compact and Hausdorff setting. 
The proof presented here originally appeared in \cite{Cas16}.

Let $X$ be a compact Hausdorff space. We already showed in Section \ref{subsection:dim.C(X)} how to prove that
$\dim_{\mathrm{nuc}}C(X) \leq \dim X$.
Before going into the proof of the opposite direction, 
let us briefly sketch it first. 

We will consider any finite open cover $\mathcal{U}$ of $X$ and a partition of unity $\left(h_r\right)_{r=1}^m$ subordinated to $\mathcal{U}$. For a sufficiently small $\epsilon$ and assuming $\dim_{\mathrm{nuc}} C(X)=n$,
we can find a finite dimensional algebra of the form $F = \bigoplus\limits_{k=0}^n F^{(k)}$, maps $\psi: C(X) \to F$ and $\varphi: F \to C(X)$, with $\varphi= \sum\limits_{k=0}^n \varphi_k$ where each $\varphi_k$ is order zero, such that they approximate the partition of unity $\left\{h_1, \ldots, h_m \right\}$.
By the structure of order zero maps (Theorem \ref{thm:structure.order.zero}), $F \cong \mathbb{C}^s$ for some $s \in \mathbb{N}$. 
We will set $f_i:=(0, \ldots,0,1,0, \ldots,0) \in \mathbb{C}^r$ where $1$ is in the $i$-th position. The family 
\begin{equation}
\mathcal{W}_0=\left\{ \varphi_k \left( f_i \right)^{-1} \left( \left( \frac{1}{m(n+1)} - \epsilon, \infty \right) \right) \mathrel{\Big|} i=1, \ldots, r; \, k=0, \ldots, n \right\}
\end{equation}
will be a cover of $X$ of order at most $n$. We will finish the proof by showing that there exists a subcover of $\mathcal{W}_0$ that refines $\mathcal{U}$.

We now do the
proof in detail.
Suppose $\dim_{\mathrm{nuc}} C(X) = n$ and let $\mathcal{U}=\left\{ U_1, \ldots, U_m \right\}$ be a finite open cover of $X$.  Without loss of generality we can assume $m \geq 2$.  Consider a partition of unity $(h_r)_{r=1}^m$ subordinated to $\mathcal{U}$.
Find $\epsilon>0$ small enough such that  
\begin{align}
\epsilon< \frac{1}{3m(n+1)}. \label{eq.sec8.initial epsilon}
\end{align}
By the choice of $\epsilon$, the following inequality holds:
\begin{align}
\frac{5\epsilon}{3} < 2\epsilon < \frac{1}{m(n+1)}- \epsilon\leq \frac{1}{m(n+1)}-\frac{\epsilon}{m(n+1)^2}. \label{eq.sec8.good estimates comm case}
\end{align}

Set $\mathfrak{F}:= \left\{1_{C(X)}, h_1, \ldots, h_m \right\}$ and using that $C(X)$ has nuclear dimension equal to $n$, we find a finite dimensional $\mathrm{C}^*$-algebra of the form $F = \bigoplus\limits_{k=0}^n F^{(k)}$, maps $\psi: C(X) \to F$ and $\varphi: F \to C(X)$ such that
\begin{enumerate}[(i)]
	\item $\psi$ is c.p.c., \label{sec8.phi.contractive}
	\item $\varphi_k := \varphi|_{F^{(k)}}$ is a c.p.c.\ order zero map,
	\item for every $x \in \mathfrak{F}$
	\begin{equation}
	\| x - \varphi \psi (x) \| < \frac{\epsilon}{m(n+1)}.
	\end{equation}
\end{enumerate}

By Theorem \ref{thm:structure.order.zero}, there exists a unital $\mathrm{C}^*$-algebra $\mathcal{M}$ containing $\varphi(F^{(k)})$ and a $^*$-homomorphism $\pi: F^{(k)} \to \mathcal{M}$ commuting with $\varphi(1_{F^{(k)}})$ such that
\begin{equation}
\varphi(a) = \varphi(1_{F^{(k)}})\pi(a) = \pi(a) \varphi(1_{F^{(k)}}).
\end{equation}
The algebra $\mathcal{M}$ given by Theorem \ref{thm:structure.order.zero} is commutative if the codomain is commutative. In fact, $\mathcal{M}$ can be assumed to be equal to $C(\beta U)$, where $U$ is an open subset of $X$ and $\beta U$ is the Stone-\v{C}ech compactification of $U$.

The algebra $F^{(k)}$ is finite dimensional, so it is a finite direct sum of matrix algebras, say $F^{(k)} = \bigoplus_{i=1}^{m_k} M_{d_i}(\mathbb{C})$.
Using that the kernel of an order zero map is an ideal of the domain, we can assume, since each $M_{d_i} \left(\mathbb{C}\right)$ is simple, that $\varphi_k|_{M_{d_i} \left(\mathbb{C}\right)}$ is injective. 
Let $a,b \in M_{d_i} \left(\mathbb{C}\right)$, using the commutativity of $\mathcal{M}$, we have
\begin{align}
\varphi(ab)  
= \varphi(1_{F^{(k)}}) \pi(a) \pi(b) 
= \varphi(1_{F^{(k)}}) \pi(b)\pi(a) 
= \varphi(ba). 
\end{align}
Since $\varphi|_{M_{d_i}(\mathbb{C})}$ is injective, we conclude that $M_{d_i}\left(\mathbb{C}\right)$ is commutative for $i=1, \ldots, m_k$. This yields $d_i=1$, and hence, $F^{(k)} = \mathbb{C}^{m_k}$ for some $m_k \in \mathbb{N}$.

In particular, for each $r\leq m$, we can write
\begin{equation}
\psi(h_r)= \left( \lambda_r^{(0)}, \ldots, \lambda_r^{(n)} \right) \in \bigoplus\limits_{k=0}^n \mathbb{C}^{m_k}
\end{equation}
where
\begin{equation}
\lambda_r^{(k)} = \left( \lambda_{1,r}^{(k)}, \ldots, \lambda_{m_k, r}^{(k)} \right) \in \mathbb{C}^{m_k}.
\end{equation}
Set $e_i^{(k)}:=(0,\ldots,0,1,0, \dots,0) \in \mathbb{C}^{m_k}$ where the $1$ is in the $i$-th position. Hence $1_{F^{(k)}}= \sum\limits_{i=1}^{m_k}e_i^{(k)}$. 
Using this notation, we can write the image of $h_r$ in the following way:
\begin{equation}\label{eq.sec8.sum.phipsi.hr}
\varphi \psi(h_r) = \varphi \left( \sum_{i=1}^{m_k} \lambda_{i,r}^{(k)} e_i^{(k)} \right)
= \sum\limits_{k=0}^{n} \sum\limits_{i=1}^{m_k} \lambda_{i,r}^{(k)} \varphi_k\left(e_i^{(k)}\right).
\end{equation}
Similarly, since $1_{C(X)}= \sum\limits_{r=1}^{m}h_r$,
\begin{align}
\varphi \psi\left(1_{C(X)}\right) = \sum\limits_{r=1}^{m}\sum\limits_{k=0}^{n} \sum\limits_{i=1}^{m_k} \lambda_{i,r}^{(k)}\varphi_k\left(e_i^{(k)}\right). \label{eq.sec8.phi psi (1) as sums of functions}
\end{align}
Observe that the functions $\varphi_k\left(e_1^{(k)}\right), \ldots, \varphi_k\left(e_{m_k}^{(k)}\right)$ are pairwise orthogonal since $\varphi_k$ is an order zero map. This implies that if $\varphi_k\left(e_i^{(k)}\right)(x) > 0$ then $\varphi\left(e_j^{(k)}\right)(x)=0$ for $j \neq i$. Hence, pointwise, the sum in (\ref{eq.sec8.phi psi (1) as sums of functions}) is at most the sum of $m(n+1)$ strictly positive summands.

For each $e_i^{(k)}$, define 
\begin{align}
W_i^{(k)} :=& \ \varphi_k\left(e_i^{(k)}\right)^{-1}\left(\left(\frac{1}{m(n+1)}- \epsilon, \infty \right)\right) \nonumber \\ 
 =& \left\{x \in X \mathrel{\Big|} \varphi_k\left(e_i^{(k)}\right) (x) > \frac{1}{m(n+1)}- \epsilon \right\}. \label{eq.sec8.def W}
\end{align}
Then set
\begin{equation}
\mathcal{W}_0: =\left\{W_i^{(k)} \mathrel{\Big|} 0 \leq k \leq n; \, 1 \leq i \leq m_k \right\}.
\end{equation}

The order of $\mathcal{W}_0$ is at most $n$ because $\varphi$ is the sum of $n+1$ order zero maps. Indeed, for each $k=0, 1, \ldots, n$, set $\mathcal{W}_0^{(0)}:=\left\{ W_i^{(0)} \mid i=1, \ldots, m_0  \right\}$ and
\begin{align}
\mathcal{W}_0^{(k)}:=&\left\{ W_i^{(k)} \mid i=1, \ldots, m_k  \right\} \setminus \bigcup_{j=0}^{k-1} \mathcal{W}_0^{(j)}, && 1 \leq k \leq n.  
\end{align} 
It is immediate that $\mathcal{W}_0 = \bigsqcup\limits_{k=0}^n \mathcal{W}^{(k)}_0$, and since the functions $\varphi_k \left( e_1^{(k)} \right), \ldots, \varphi_k \left( e_{m_k}^{(k)} \right)$ are pairwise orthogonal, we have that $W_i^{(k)} \cap W_j^{(k)} = \emptyset$ for $i \neq j$.

Momentarily fix $W_j^{(s)}$ for some $j$ and $s$. We will show that if it is not contained in some element $U_r$ of the original cover $\mathcal{U}$, then the corresponding coefficient $\lambda_{j,r}^{(s)}$ in equation \eqref{eq.sec8.phi psi (1) as sums of functions} is ``small''. 
Precisely, suppose there exists $U_r \in \mathcal{U}$ such that $W_j^{(s)}\cap\left(X \setminus U_r \right)\neq \emptyset$, we will prove  $\lambda_{j,r}^{(s)} < 5 \epsilon / 3$. 

Let $x \in W_j^{(s)}\cap\left(X \setminus U_r \right)$. Using that the support of $h_r$ is contained in $U_r$, we obtain $h_r(x)=0$. 
By hypothesis
\begin{align}
\left\|h_r - \varphi \psi (h_r)\right\| < \frac{\epsilon}{m(n+1)}, \label{eq.sec8.commutative case eq 1}
\end{align} 
then
\begin{align}
\sum\limits_{k=0}^{n} \sum\limits_{i=1}^{m_k} \lambda_{i,r}^{(k)} \varphi_k\left(e_i^{(k)}\right) (x) \overset{\eqref{eq.sec8.sum.phipsi.hr}}{=} \varphi \psi(h_r)(x) < \frac{\epsilon}{m(n+1)}.
\end{align}
Observe that all summands are positive, hence
\begin{align}
\lambda_{j,r}^{(s)} \varphi_s\left(e_j^{(s)}\right)(x) < \frac{\epsilon}{m(n+1)}. \label{eq.sec8.lambda phi  < epsilon/m(n+1)}
\end{align}
Using that $x \in W_j^{(s)}$ and
\begin{align}
\epsilon \overset{(\ref{eq.sec8.initial epsilon})}{<} \frac{1}{3m(n+1)} < \frac{2}{5m(n+1)}, \label{eq.sec8.epsilon < 1/3<2/5}
\end{align}
we obtain
\begin{align}
\frac{\epsilon}{m(n+1)} & \overset{(\ref{eq.sec8.lambda phi  < epsilon/m(n+1)})}{>} \lambda_{j,r}^{(s)} \varphi_s\left(e_j^{(s)}\right)(x) \overset{(\ref{eq.sec8.def W})}{>} \lambda_{j,r}^{(s)}\left(\frac{1}{m(n+1)}-\epsilon\right) \nonumber \\
& \overset{\eqref{eq.sec8.epsilon < 1/3<2/5}}{>} \lambda_{j,r}^{(s)}\left(\frac{1}{m(n+1)}- \frac{2}{5m(n+1)} \right) = \frac{3\lambda_{j,r}^{(s)}}{5m(n+1)}. 
\end{align}
Thus
\begin{align}
\lambda_{j,r}^{(s)} < \frac{5\epsilon}{3}. \label{eq.sec8.lambda_ir is small enough}
\end{align}

Next step is proving that for any $x \in X$
there is at least one coefficient $\lambda_{j,r}^{(s)}$ which is ``large enough'', i.e.\ $\lambda_{j,r}^{(s)} \geq 5\epsilon /3$. 

Let $x \in X$ and $h_r$ such that $h_r(x) \geq \frac{1}{m}$. Such function exists by the \textit{pigeonhole principle} and because $(h_1, \ldots, h_m)$ is a partition of unity.
Then, by (\ref{eq.sec8.commutative case eq 1}), we have
\begin{align}
\sum\limits_{k=0}^{n} \sum\limits_{i=1}^{m_k} \lambda_{i,r}^{(k)} \varphi_k\left(e_i^{(k)}\right) (x) \overset{\eqref{eq.sec8.sum.phipsi.hr}}{=} \varphi \psi(h_r)(x) \geq \frac{1}{m}-\frac{\epsilon}{m(n+1)}. 
\end{align}
Then, by the pigeonhole principle again, there exists $s \in \{0, \ldots, n \}$ such that
\begin{align}
\sum\limits_{i=1}^{m_s} \lambda_{i,r}^{(s)} \varphi_s\left(e_i^{(s)}\right) (x) = \varphi_s \psi(h_r)(x) \geq \frac{1}{m(n+1)}-\frac{\epsilon}{m(n+1)^2}. 
\end{align}

The functions $\varphi_s\left(e_1^{(s)}\right), \ldots, \varphi_s\left(e_{m_s}^{(s)}\right)$ are pairwise orthogonal, thus all the summands in 
the previous inequality are equal to zero but one. Hence, there exists $j \in\{1, \ldots, m_s \}$ such that
\begin{align}
\lambda_{j,r}^{(s)} \varphi_s\left(e_j^{(s)}\right) (x) =\sum\limits_{i=1}^{m_s} \lambda_{i,r}^{(s)} \varphi_s\left(e_i^{(s)}\right) (x) \geq \frac{1}{m(n+1)} - \frac{\epsilon}{m(n+1)^2}. 
\end{align}
Since $\lambda_{j,r}^{(s)}$ and $\varphi_s\left(e_j^{(s)}\right) (x)$ are positive numbers less than $1$, we obtain
\begin{align}
\lambda_{j,r}^{(s)} \geq \frac{1}{m(n+1)} - \frac{\epsilon}{m(n+1)^2} 
\label{eq.sec8.good estimates comm case II}
\end{align}
and
\begin{align}
\varphi_s\left(e_j^{(s)}\right) (x) \geq \frac{1}{m(n+1)} - \frac{\epsilon}{m(n+1)^2} > \frac{1}{m(n+1)} - \epsilon.
\end{align}
This immediately shows $x \in W_{j}^{(s)}$.  
Moreover, we have
\begin{equation}
\lambda_{j,r}^{(s)} \overset{\eqref{eq.sec8.good estimates comm case II}}{>} \frac{1}{m(n+1)} - \frac{\epsilon}{m(n+1)^2}\overset{\eqref{eq.sec8.good estimates comm case}}{>}  
\frac{5 \epsilon}{3}.
\end{equation}
Then, by (\ref{eq.sec8.lambda_ir is small enough}), 
$W_{j}^{(s)} \cap \left(X \setminus U_r \right) = \emptyset$. In other words, $W_{j}^{(s)} \subseteq U_r$. 

This shows that for any $x \in X$, there exists some $U_r$ and $W_{j}^{(s)}$ such that $x \in W_{j}^{(s)} \subseteq U_r$. Set 
\begin{align}
\mathcal{W}= \left\{ W_{i}^{(k)} \in \mathcal{W}_0 \; \middle| \; W_{i}^{(k)} \subseteq U_r \; \mbox{for some }r=1, \ldots,m \right\}.
\end{align}
Our previous arguments show that $\mathcal{W}$ is indeed a cover of $X$ that is contained in $\mathcal{W}_0$. Hence, its order is at most $n$ and, by construction, $\mathcal{W}$ refines $\mathcal{U}$. Therefore 
\begin{align}
\dim X \leq \dim_{\mathrm{nuc}} C(X). 
\end{align}

The argument from Section \ref{subsection:dim.C(X)} actually shows that $\mathrm{dr} \, C(X) \leq \dim X$.
Finally, using that nuclear dimension is always smaller than decomposition rank yields the following,
\begin{equation}
\mathrm{dr} \, C(X) \leq \dim X \leq \dim_{\mathrm{nuc}} C(X) \leq \mathrm{dr} \, C(X)
\end{equation} 

If we want to include the second countable locally compact case as well, we need to use the following identities regarding a locally compact space $X$ and its one point compactification $\alpha X$ (\cite[Remark 2.11]{WZ10} \cite[Corollary 3.5.8]{Pe75}),
\begin{align}
\dim_{\mathrm{nuc}} C_0(X)  =  \dim_{\mathrm{nuc}} C(\alpha X), \qquad \text{ and } \qquad \dim X = \dim \alpha X.
\end{align}
We then obtain the following theorem.

\begin{theorem}
	Let $X$ be a second countably locally compact Hausdorff space. Then
	\begin{equation}
	\dim X = \dim_{\mathrm{nuc}} C_0(X) = \mathrm{dr} \, C_0(X).
	\end{equation}
\end{theorem}


\newcommand{\etalchar}[1]{$^{#1}$}

\end{document}